\documentclass[a4paper,reqno]{amsart}

\usepackage{amsthm,amssymb,amsmath,amsfonts}
\usepackage{verbatim}
\usepackage{bbm}
\usepackage{hyperref}
\usepackage[all]{xy}
\usepackage{latexsym}
\usepackage{color}
\usepackage{tikz, tikz-qtree}
\usepackage{graphicx}
\usepackage{caption, subcaption}
\usepackage{verbatim}
\usepackage{rotating} 
\usepackage[colorinlistoftodos]{todonotes} 

\usepackage{amsaddr}

\usepackage{natbib}

\theoremstyle{definition}

\def\R{\mathbb{R}}

\begin{document}

\title{An Interpolating Family of size distributions}
\author{Corinne Sinner$^*$,  Yves Dominicy$^\dagger$,  Christophe Ley$^\ddagger$, Julien~Trufin$^*$ and Patrick Weber$^*$}
\address{$^*$Department of Mathematics, Universit\'e libre de Bruxelles (ULB), Belgium, E-Mail: csinner@ulb.ac.be, julien.trufin@ulb.ac.be, pweber@ulb.ac.be.}
\address{$^\dagger$SBSEM, ECARES, Universit\'e libre de Bruxelles (ULB), Belgium, E-Mail:~yves.dominicy@ulb.ac.be.}
\address{$^\ddagger$Department of Applied Mathematics, Computer Science and Statistics, Ghent~University, Belgium, \\
E-Mail: christophe.ley@ugent.be.}
%

\maketitle

\begin{abstract}
We introduce a new five-parameter family of size distributions on the semi-finite interval $[x_0,\infty), x_0\geqslant0,$ with two attractive features. First, it interpolates between power laws, such as the Pareto distribution, and power laws with exponential cut-off, such as the 
Weibull distribution. The proposed family  is thus very flexible and spans over a broad range of  well-known size distributions which are special cases of our family. Second, it has important tractability advantages over the popular five-parameter Generalized Beta distribution. We derive the hazard function, survival function, modes and quantiles, propose a random number generation procedure and discuss  maximum likelihood estimation issues. Finally, we illustrate the wide applicability and fitting capacities of our new model on basis of three  real data sets  from very diverse domains, namely actuarial science, environmental science and survival analysis. 
\end{abstract}

\hfill

\noindent%


\section{Introduction}
\label{sec:intro}

Size distributions are probability laws designed to model data that only take positive values. Typical examples of such size-type data are claim sizes in actuarial science or wind speeds in meteorology. Nonetheless, the spectrum of application areas is much broader and positive observations appear naturally in survival analysis~\citep{lawless2003, leewang2003}, environmental science~\citep{marchenko2010}, network traffic modeling~\citep{mitzenmacher2004}, reliability theory~\citep{rausandhoyland2004}, economics \citep{eeckhout2004, luttmer2007, gabaix2016} and hydrology~\citep{clarke2002}. Parametric distributions are a simple and effective way to convey the information contained in those data. Given the range of distinct domains of application, it is not surprising that there exists a plethora of different size distributions and that the quest for the right size distribution in a given situation has developed into a very active research area over the past years~\citep{ortega2015, hongrubio2016, asgharzadeh2016}.

\hfill

What desirable properties should a  size distribution possess? Obviously, it should be quite flexible, meaning that it is able to model very diverse data shapes, yet it should ideally remain of a tractable form. Moreover, it should be parsimonious in terms of the parameters it uses, and each parameter should bear a clear interpretation. The latter points, tractability and interpretability, are especially important to practitioners from other domains.

\hfill

A very popular size distribution is the Pareto distribution, also called Pareto type I distribution, with probability density function 
$$x \mapsto  \frac{\alpha x_{\mathrm{0}}^{\alpha}}{x^{\alpha+1}}, \quad x \in [x_0,\infty),$$
where $x_0\geqslant0$ is a location parameter and $\alpha>0$ is a shape parameter known as the tail or Pareto index. Vilfredo Pareto used this law to model the distribution of income as well as the allocation of wealth among individuals~\citep{pareto1964}. Over the years, the Pareto law has further been applied to city sizes, file size distribution of internet traffic, sizes of meteorites, or the sizes of sand particles~\citep{reed2004}.

\hfill

The Pareto distribution is a member of the power laws, which are typically of the form $x\mapsto k x^{-\alpha}$, with normalizing constant $k$ and power exponent $\alpha >0$. 
Power law distributions are employed in a vast  range of situations, such as the modeling of the number of hits on web pages~\citep{huberman1999}, income of top earners in areas of arts, sports and business~\citep{rosen1981}, and inequalities of income and wealth~\citep{piketty2014, toda2015}. For further information and references, we refer the interested reader to~\cite{sornette2003} and~\cite{mitzenmacher2004}. 

\hfill

A popular alternative to power laws are power laws with exponential cut-off,  whose densities take the form $x \mapsto k x^{-\alpha}\mathrm{e}^{-\beta x}$ with normalizing constant $k$, power exponent $\alpha >0$ and rate parameter $\beta>0$.
A power law with exponential cut-off behaves like a power law for small values of $x$, while its tail behavior is governed by a decreasing exponential. A famous representative of this class of distributions is the Weibull distribution with density
$$x \mapsto \frac{\alpha}{\sigma}\left(\frac{x}{\sigma}\right)^{\alpha-1}e^{-\left(\frac{x}{\sigma}\right)^{\alpha}}, \quad x \in [0,\infty),$$
where $\alpha > 0$ is a shape parameter regulating tail-weight and $\sigma > 0$ is a scale parameter. 
While the Weibull distribution already appeared in~\cite{rosin1933} to describe the distribution of particle sizes, it gained its prominence and name after Waloddi Weibull who showed in 1951 that the distribution can be successfully applied to seven very different case studies. Nowadays, the Weibull distribution is widely used in various domains such as life data analysis~\citep{nelson2005}, wind speed modeling~\citep{manwell2009} and hydrology~\citep{clarke2002}.

\hfill

\begin{figure}[h!]
\begin{center}
\includegraphics[scale=0.5]{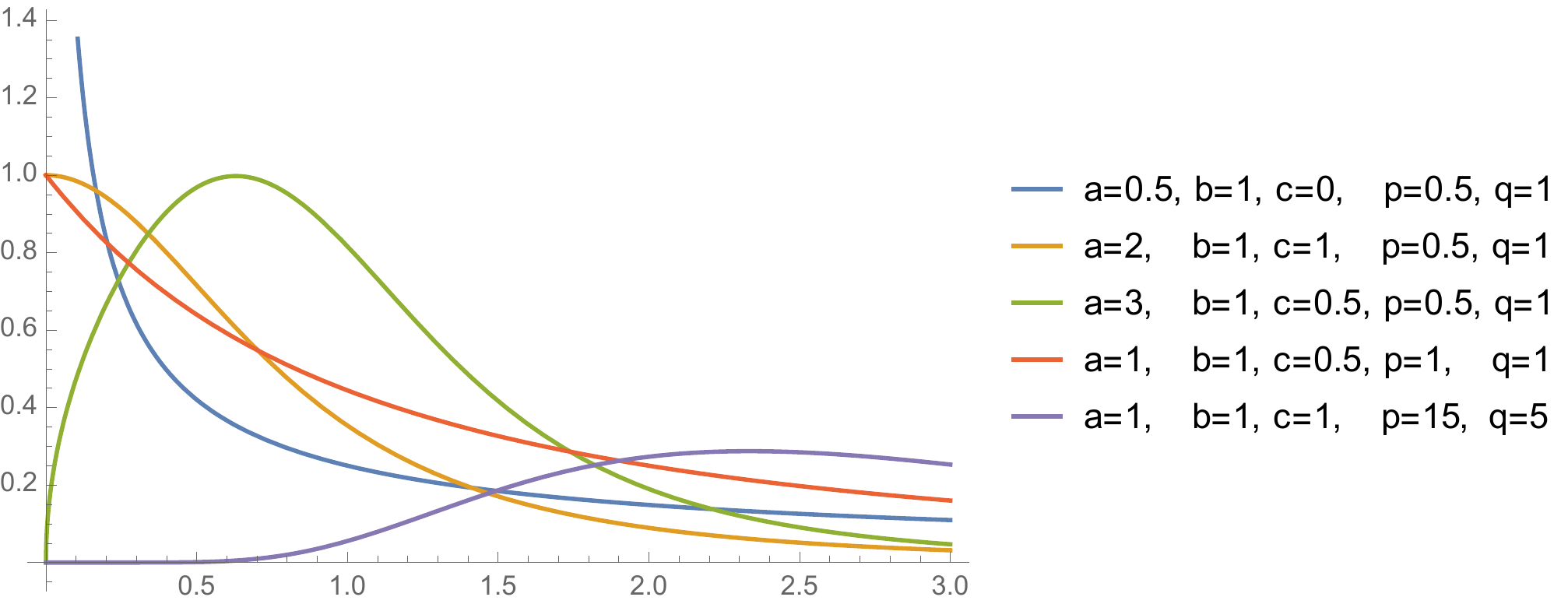}
\end{center}
\caption{Density plots of the GB distribution for various choices of parameters.} \label{GBplot}
\end{figure}

Given the wealth of different size distributions, the practitioner is often confronted to the question ``Which distribution shall I use in what situation?". The disparity between the Pareto and the Weibull, and more generally between power laws and power laws with exponential cut-off, renders this choice even more delicate. A solution to this dilemma has been proposed by~\cite{mcdonaldxu1995} who introduced the Generalized Beta (GB) distribution with density 
$$\text{GB}(x;a,b,c,p,q) = \frac{|a|x^{ap-1}(1-(1-c)(\frac{x}{b})^{a})^{q-1}}{b^{ap}B(p,q)(1+c(\frac{x}{b})^{a})^{p+q}}, \quad x^{a} \in \left( 0,\frac{b^a}{1-c} \right).$$
Here, $B(p,q)$ stands for the beta function, the parameter $c$ satisfies $0 \leqslant c \leqslant 1$, the shape parameter $a$ is non-zero whilst the scale parameter $b$, the shape parameter $p$ as well as the skewness parameter $q$ are all positive. The GB is extremely flexible and incorporates more than thirty known distributions as special cases, among which the Pareto and Weibull. In particular, it encompasses the already broad four-parameter subfamilies GB1 and GB2 suggested by~\cite{mcdonald1984}. These respectively correspond to $\text{GB1}(x;a,b,p,q) = \text{GB}(x;a,b,c=0,p,q)$ and $\text{GB2}(x;a,b,p,q) = \text{GB}(x;a,b,c=1,p,q)$. The GB distribution is attractive because of this versatility, as it enables different shapes (see Figure \ref{GBplot}) and hence avoids the researcher or practitioner to try out various choices of distributions. 

However, this family of distributions also suffers from serious flaws such as the absence of a  tractable cumulative distribution function, entailing a very intricate quantile structure and rendering random variable generation difficult. For further information about the GB and a general overview on size distributions, we refer the interested reader to~\cite{mcdonaldxu1995} and~\cite{kleiberkotz2003}.
 
\hfill

In the present paper, we propose a viable alternative to the very popular GB distribution that shares many of its advantages  but overcomes its major shortcomings. The new size distribution has the density
\begin{equation}\label{newdist}
	\text{IF}(x;p,b,c,q,x_0)=\frac{|b|q}{c}\left(\frac{x-x_0}{c}\right)^{b-1} G_p(x)^{-q-1}\left(1- \frac{1}{p+1}G_p(x)^{-q}\right)^p, 
\end{equation}
for $x\in [x_0,\infty)$, and with $p\in[0,\infty], b\neq0, c,q>0, x_0\geqslant 0$ and
\begin{equation*}
	G_p(x)=(p+1)^{-\frac{1}{q}}+ \left(\frac{x-x_0}{c}\right)^b.
\end{equation*}
The roles of the various parameters at play will be laid out in Section~\ref{sec:para}. The alerted reader will certainly already have recognized the densities of various size distributions in~\eqref{newdist}, including the Pareto and Weibull. The originality of our proposal rests on its construction: we are using a technique from statistical mechanics (see Section 2) that allows us to interpolate between the Pareto and Weibull distributions, even more generally, between power laws and power laws with exponential cut-off. Thus, we are precisely  finding a path from one end of the spectrum to the other, and this moreover in a constructive way. 

\hfill

The remainder of the paper is organized as follows. In Section 2 we describe the new five-parameter family of size distributions, explain how it interpolates between power laws and power laws with exponential cut-off and elucidate the role of each of the five parameters. In Section 3 we summarize some of the special cases belonging to the Interpolating Family of size distributions. Section 4 presents the main properties of the family such as the cumulative distribution function, survival function, quantile function, random number generation procedure, moments and mode, while Section 5 contains a   discussion on  inferential aspects. In Section 6 we analyze three  real data sets from distinct areas, namely  actuarial science, environmental science and survival analysis, to show the flexibility and effectiveness of the newly introduced family of size distributions. Section 7 concludes and technical derivations are presented in the Appendix.


\section{The Interpolating Family: construction and parameter interpretation} 
\label{sec:IF}

In this section we present in detail our new five-parameter family of size distributions, which we call \textit{Interpolating Family} for reasons that will become obvious through the construction described in Section~\ref{sec:meth}. Section~\ref{sec:para} expounds on the role of each of the five parameters. 


\subsection{Construction of the family} 
\label{sec:meth}

As announced in the Introduction, our goal is to build a size distribution which incorporates both power laws and power laws with exponential cut-off. 
To show that~\eqref{newdist} indeed satisfies this requirement, we start by writing up power law distributions and  power law distributions with exponential cut-off in a unified language.

\subsubsection{Power laws}

The probability density function (pdf) of a typical power law distribution corresponds to
$$x\mapsto q(1+x)^{-q-1}, \quad x \in [0,\infty),$$
where the tail behavior is governed by the shape parameter $q>0$. To get a more flexible distribution, one may add various parameters, such as a scale parameter $c>0$, a location parameter $x_0 \geqslant 0$ and/or a shape parameter $b>0$, leading to
\begin{equation*}
	x\mapsto \frac{bq}{c} \left(\frac{x-x_0}{c}\right)^{b-1} \left(1+\left(\frac{x-x_0}{c}\right)^b\right)^{-q-1}, \quad x \in [x_0,\infty).
\end{equation*}
Alternatively, in terms of the function $G_0(x)=1+\left(\frac{x-x_0}{c}\right)^b$, the pdf can be written under the form
\begin{equation*}
	f_0(x)=q~g_0(x) G_0(x)^{-q-1},  \quad x \in [x_0,\infty), 
\end{equation*}
where $g_0(x)=\frac{\mathrm{d}}{\mathrm{d}x}G_0(x)=\frac{b}{c}\left(\frac{x-x_0}{c}\right)^{b-1}$. We point out that the function $G_0(x)$ has been chosen such that the following boundary conditions are satisfied: $G_0(x_0)=1$ and ${\lim \limits_{x \to \infty} G_0(x)=\infty}$.

\subsubsection{Power laws with exponential cut-off}

The pdf of a typical power law distribution with exponential cut-off reads
\begin{equation*}
	x\mapsto qx^{-q-1}e^{-x^{-q}}, \quad x \in [0,\infty).
\end{equation*}
The shape parameter $q>0$ still controls the tail behavior and, just as for power laws, we may increase the flexibility of the model by adding scale, location and shape parameters to get
\begin{equation*}
	x\mapsto\frac{bq}{c} \left(\frac{x-x_0}{c}\right)^{-bq-1} e^{-\left(\frac{x-x_0}{c}\right)^{-bq}}, \quad x \in [x_0,\infty).
\end{equation*}
Alternatively, we may write the pdf in terms of the function $G_{\infty}(x)=\left(\frac{x-x_0}{c}\right)^b$ as
\begin{equation*}
	f_{\infty}(x)=q~g_{\infty}(x) G_{\infty}(x)^{-q-1} e^{-G_{\infty}(x)^{-q}}, \quad x \in [x_0,\infty), 
\end{equation*}
where $g_{\infty}(x)=\frac{\mathrm{d}}{\mathrm{d}x}G_{\infty}(x)=\frac{b}{c}\left(\frac{x-x_0}{c}\right)^{b-1}$. Note that the function $G_{\infty}(x)$ has been chosen such that ${G_{\infty}(x_0)=0}$ and $\lim \limits_{x \to \infty} G_{\infty}(x)=\infty$.

\subsubsection{Interpolating Family} 

If we want a highly flexible distribution including both power laws and power laws with exponential cut-off, we need a way to build densities interpolating between  $f_0(x)$ and $f_{\infty}(x)$. To  this end, we introduce a mild variant of the one-parameter deformation of the exponential function popularized in the seminal paper \cite{tsallis1988} in the context of non-extensive statistical mechanics. A more detailed account can be found in  the review paper~\cite{tsallis2002}.

\hfill

For any $p \in [0,\infty]$, we define the $p$-exponential\footnote{We would like to point out that the classical $q$-exponential defined by Tsallis has the form $\tilde{e}_q(x)= \left(1+(1-q)x\right)^{\frac{1}{1-q}}$. 
We slightly modified the deformation path in order to simplify the calculations.} by
\begin{equation*}
	e_p(x)=\left(1-\frac{1}{p+1} x\right)^p, \quad x \in [0,p+1].
\end{equation*}
The extreme cases $p=0$ and  $p \to  \infty$ respectively correspond to  $e_0(x)=1$ over~$[0,1]$ and $e_\infty(x)=e^{-x}$ over $[0,\infty)$. With this in mind, it is natural to consider densities of the type
\begin{equation}\label{firsttry}
	f_p(x)= q~g_p(x) G_p(x)^{-q-1} e_p\left(G_p(x)^{-q}\right), \quad x \in [x_0,\infty), 
\end{equation}
with $g_p(x)=\frac{\mathrm{d}}{\mathrm{d}x}G_p(x)$, where we have not defined the function $G_p(x)$ yet. Just as~$e_p(x)$ interpolates between 1 and $e^{-x}$, the mapping $G_p$ should also vary between~$G_0$ and $G_\infty$. 
Hence, with the parameters $c>0$, $b>0$ and $x_0 \geqslant 0$ bearing the same interpretation as before, the map $G_p$ could be chosen as $G_p(x)=k+\left(\frac{x-x_0}{c}\right)^b$ for some constant $k$.
A quick calculation shows that $k=(p+1)^{-\frac{1}{q}}$ is the right choice for $f_p(x)$ to integrate to one over its domain. Consequently
\begin{equation*}
	G_p(x)=(p+1)^{-\frac{1}{q}}+ \left(\frac{x-x_0}{c}\right)^b, \quad x \in [x_0, \infty).
\end{equation*}
Since $G_p(x)^{-q}$ with $q>0$  maps $[x_0, \infty)$ onto $[0, p+1]$, the function $e_p(G_p(x)^{-q})$ is well-defined. The pointwise convergence of the resulting density $f_p$ to $f_0$ as $p$ tends to zero (respectively $f_p$ to $f_\infty$ as $p\to \infty$) can be shown by straightforward limit calculations which we omit here. 

\hfill

The density~\eqref{firsttry} now almost corresponds to the density announced in the Introduction. Relaxing the condition $b>0$ into $b\in\R_0$, we finally end up with 
\begin{equation}\label{IF}
\text{IF}(x;p,b,c,q,x_0)=\text{sign}(b)q~g_p(x) G_p(x)^{-q-1} e_p\left(G_p(x)^{-q}\right), \quad x\in [x_0,\infty).
\end{equation}
The relaxation on $b$ only entails a minor change in the normalizing constant, which remains extremely simple. We call IF the \emph{Interpolating Family} of size distributions. The density depends on five parameters $p,b,c,q$ and $x_0$, which we will discuss in more detail in the next section. 


\subsection{Interpretation of the parameters}\label{sec:para} 

\begin{figure}[h!]
	\centering
	\begin{subfigure}{0.48\textwidth}
		\includegraphics[scale=0.5]{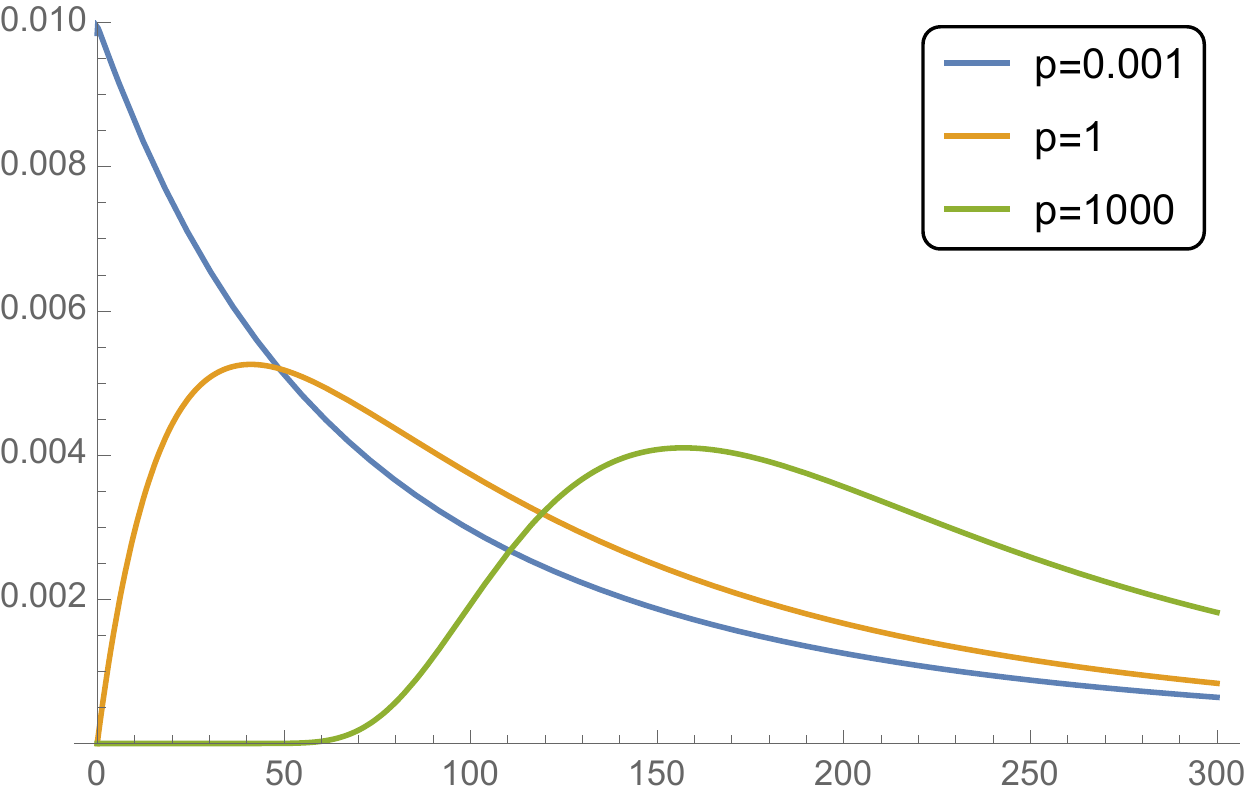}
	\end{subfigure}
	\hfill
	\begin{subfigure}{0.48\textwidth}
		\includegraphics[scale=0.5]{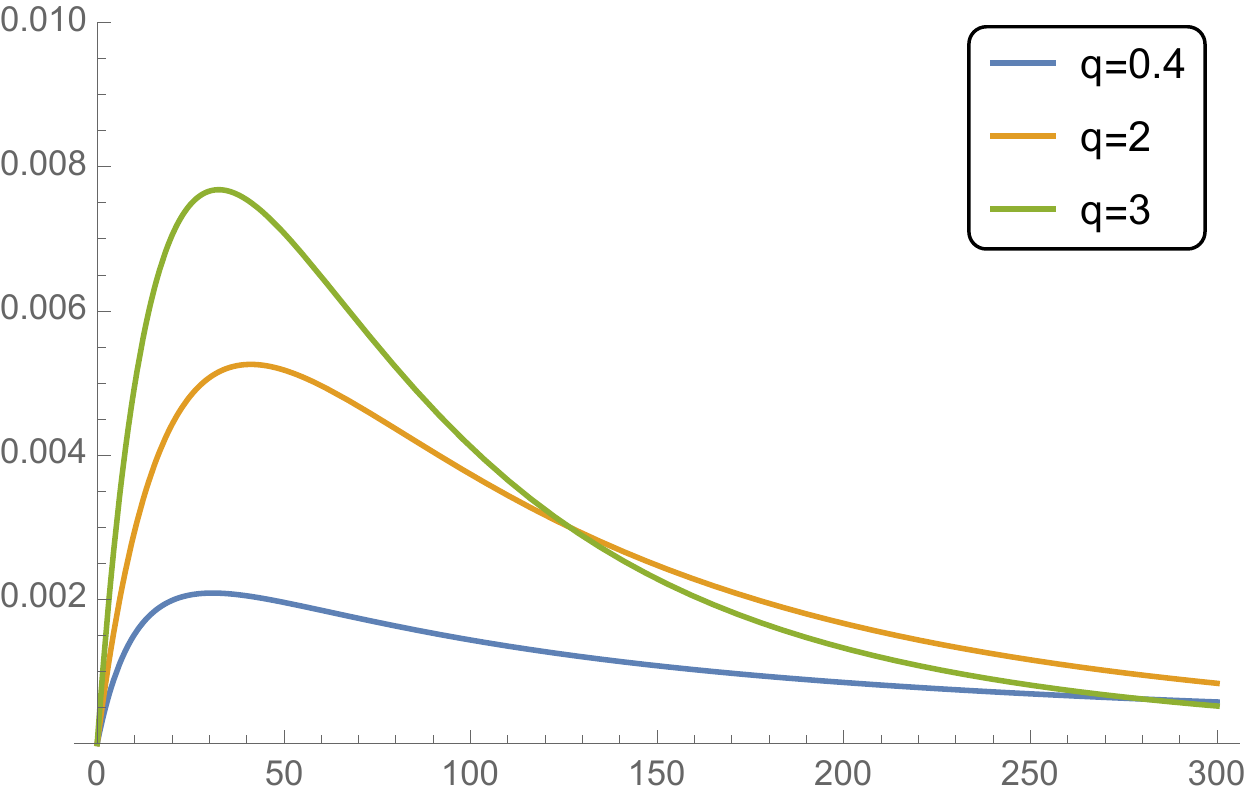}
	\end{subfigure}
	\centering
	\begin{subfigure}{0.48\textwidth}
		\includegraphics[scale=0.5]{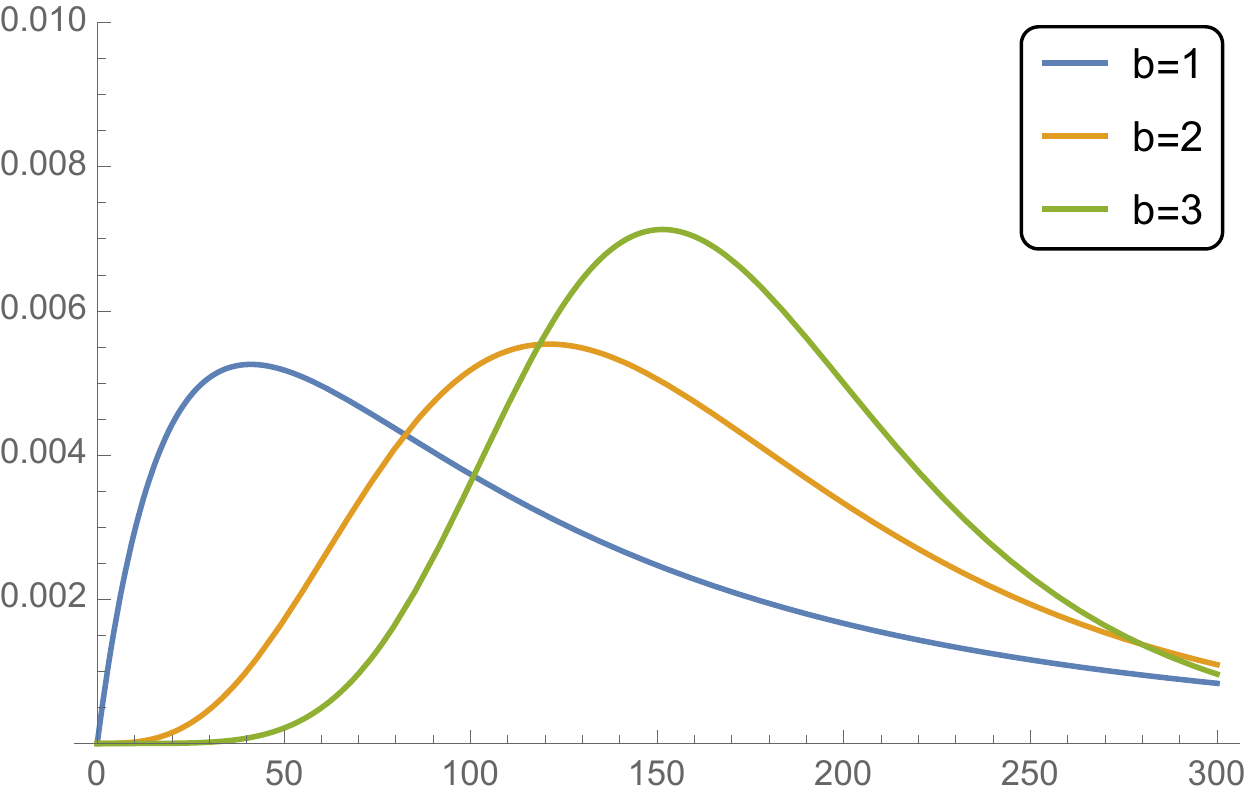}
	\end{subfigure}
	\hfill
	\begin{subfigure}{0.48\textwidth}
		\includegraphics[scale=0.5]{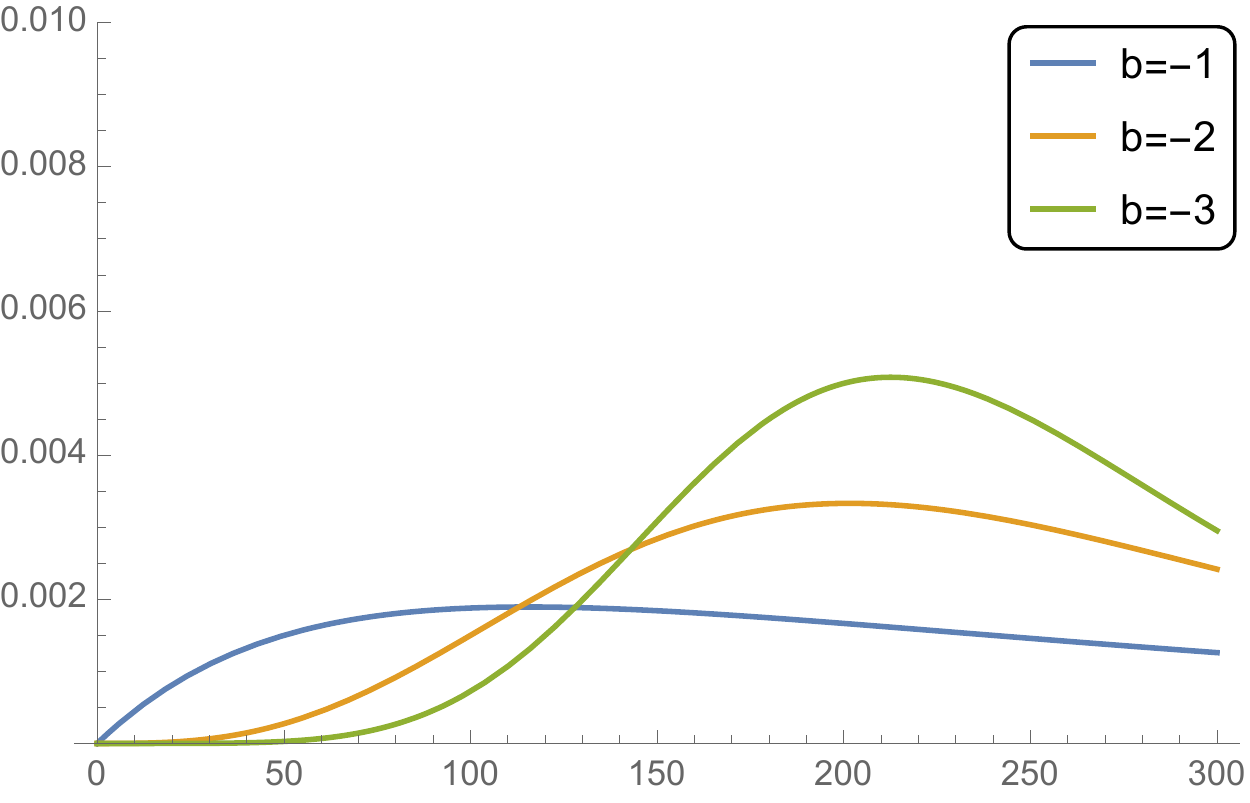}
	\end{subfigure}
	\centering
	\begin{subfigure}{0.48\textwidth}
		\includegraphics[scale=0.5]{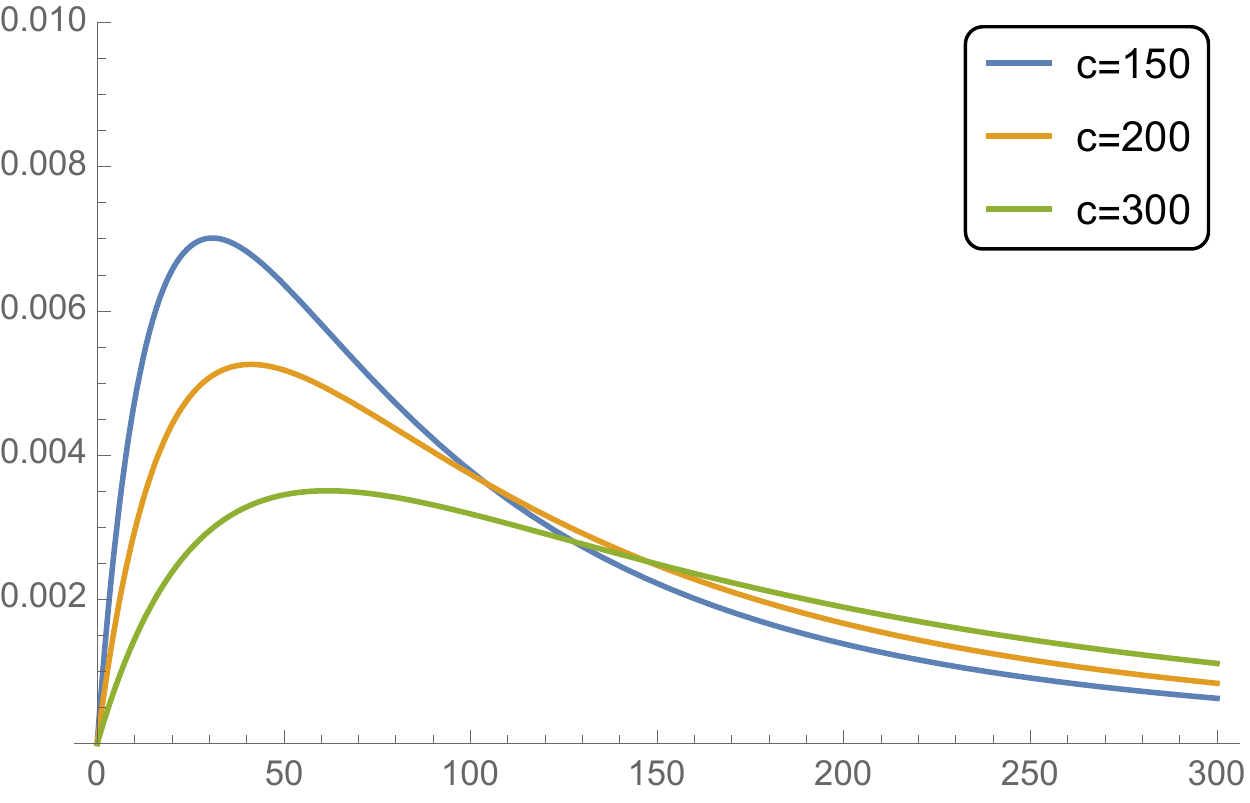}
	\end{subfigure}
	\hfill
	\begin{subfigure}{0.48\textwidth}
		\includegraphics[scale=0.5]{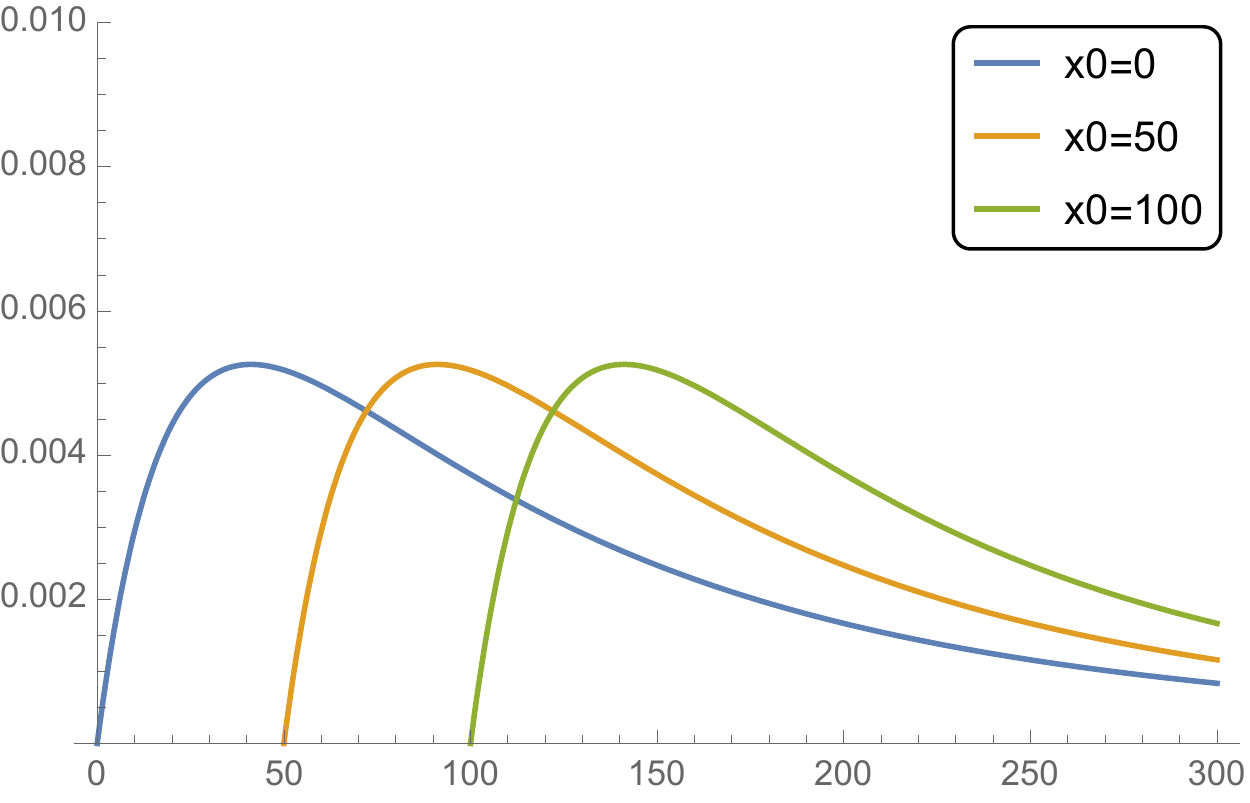}
	\end{subfigure}
	\caption{ \label{Parameters} Density plots of the IF distribution where, except for the parameter we vary in each plot, the parameters remain fixed to $p=1$, $b=1$, $c=200$, $q=2$ and $x_0=0$.}
\end{figure}

For the sake of illustration, we provide density plots of the IF distribution in Figure \ref{Parameters}. Except for the parameter we are varying, all the parameters remain fixed to $p=1$, $b=1$, $c=200$, $q=2$ and $x_0=0$. 

\newpage

Figure~\ref{Parameters} provides a visual inspection of the roles endorsed by the five parameters: $x_0\geqslant0$ is a location parameter (smaller than the lowest value of the data), $c>0$ a scale parameter, $q>0$ a tail-weight parameter, and $b\in\R_0$ a shape parameter regulating the skewness. By changing the sign of $b$ in IF, one gets  the  Inverse-IF distribution (such as, for instance, the Rayleigh and Inverse Rayleigh distribution, see below). A crucial role is played by $p\in[0,\infty]$ as it enables us to interpolate between power laws and power laws with exponential cut-off. We therefore name it \textsl{interpolation parameter}.


\section{Special cases and three main IF subfamilies} \label{sec:special}

One major appeal of the IF distribution is that it contains a plethora of well-known size distributions as special cases. For a clearer structure, we define three four-parameter subfamilies: 
\begin{itemize}
\item{the IF1 distribution where $p=0$,}
\item{the IF2 distribution where $p \to \infty$,}
\item{the IF3 distribution where $p\in(0,\infty)$ and $b=1$.}
\end{itemize}

Of course, there remain several other parameter combinations in the Interpolating Family, and perhaps in the future other interesting subfamilies will be given special attention.


\subsection*{The IF1 distribution}

In the power law limit $p=0$, the pdf of the resulting four-parameter family of distributions, called \emph{Interpolating Family of the first kind~(IF1)}, is given by
$$
	f_0(x)=\text{sign}(b)q~g_0(x)G_0(x)^{-q-1}		=\frac{|b|q}{c} \left(\frac{x-x_0}{c}\right)^{b-1} \left(1+\left(\frac{x-x_0}{c}\right)^{b}\right)^{-q-1},
$$
where $x\in [x_0,\infty)$. Special cases of the IF1 distribution are, in decreasing order of the number of parameters,
the Lindsay--Burr type III distribution ($b<0$), the Pareto type IV distribution ($b>0$), the Dagum distribution ($b<0$ and~$x_0=0$), the Pareto type II distribution ($b=1$), the Pareto type III distribution ($b>0$ and~$q=1$), the Tadikamalla--Burr type XII distribution ($b>0$ and~$x_0=0$), the Pareto type I distribution ($b=1$ and~$c= x_0>0$), the Lomax distribution ($b=1$ and~$x_0=0$), the Burr type XII distribution ($b>0, c=1$ and~$x_0=0$) and the Fisk distribution ($b>0, q=1$ and~$x_0=0$).

\subsection*{The IF2 distribution}

In the power law with exponential cut-off limit $p \to \infty$, the pdf of the resulting four-parameter family of distributions, called \emph{Interpolating Family of the second kind (IF2)}, is given by
$$f_{\infty}(x)= \text{sign}(b)q~g_{\infty}(x)G_{\infty}(x)^{-q-1} e^{-G_{\infty}(x)^{-q}}
			= \frac{|b|q}{c} \left(\frac{x-x_0}{c}\right)^{-bq-1} e^{-\left(\frac{x-x_0}{c}\right)^{-bq}},
$$
where $x\in[x_0,\infty)$. Special cases of the IF2 distribution are, in decreasing order of the number of parameters, the Weibull distribution ($b=-1$; if moreover $x_0=0$, we find the two-parameter Weibull distribution), the Fr\'{e}chet distribution ($b=1$; if moreover $x_0=0$, we find the two-parameter Fr\'{e}chet distribution), the Gumbel type~II distribution ($b=1$ and $x_0=0$), the Rayleigh distribution ($b=-1,q=2$ and $x_0=0$), the Inverse Rayleigh distribution ($b=1,q=2$ and $x_0=0$), the Exponential distribution ($b=-1,q=1$ and $x_0=0$), and the Inverse Exponential distribution ($b=1,q=1$ and $x_0=0$).


\subsection*{The IF3 distribution}

The \emph{Interpolating Family of the third kind (IF3)} is characterized by  $0<p<\infty$ and $b=1$, resulting in the pdf
$$f_{p,1}(x) = \frac{q}{c} \left((p+1)^{-\frac{1}{q}} +\frac{x-x_0}{c}\right)^{-q-1} \left(1-\frac{1}{p+1}\left((p+1)^{-\frac{1}{q}} + \frac{x-x_0}{c}\right)^{-q}\right)^p,$$
where $x\in [x_0,\infty)$. 
Special cases   of the IF3 distribution are  the Generalized Lomax distribution ($x_0=0$) and the Stoppa distribution ($x_0=c(p+1)^{-\frac{1}{q}}$).


\subsection*{Distribution tree} 

A visual summary of the structure inherent to the Interpolating Family of distributions with its various special cases is given in the \emph{distribution tree} depicted below. Since the inverse of each distribution is obtained by switching the sign of the parameter $b$, we only give the tree for positive values of $b$. 

\begin{center}
\begin{sideways}
\begin{tikzpicture}
 \tikzstyle{block}=[rectangle, draw=blue!40, thick, fill=blue!10, text width=6em, text centered, rounded corners, minimum height=2em];
 \tikzstyle{indic}=[text width=6em];
 \tikzstyle{indic2}=[text width=15em];  
  \path (1.75,0) node [block] (IF) {IF}
            (-4.75,-3) node [block] (IF1) {IF1 \quad (Pareto IV)}
	   (1.75,-3) node [block] (F) {IF3}
            (8.25,-3) node [block] (IF2) {IF2}
            (-8,-6) node [block] (Par3) {Pareto III}
            (-4.75,-6) node [block] (TB12) {T-Burr XII}
            (-1.5,-6) node [block] (Par2) {Pareto II}
            (1.75,-6) node [block] (GL) {Generalized Lomax}
            (5,-6) node [block] (Stoppa) {Stoppa}
            (8.25,-6) node [block] (Frechet) {Fr\'{e}chet}
            (-8,-9) node [block] (Fisk) {Fisk}
            (-4.75,-9) node [block] (Burr3) {Burr XII}
            (-1.5,-9) node [block] (Lomax) {Lomax}
            (2.5,-9) node [block] (Par1) {Pareto I}
            (6.62,-9) node [block] (Gumbel) {Gumbel II}
            (5,-12) node [block] (IExponential) {Inverse Exponential}
            (8.25,-12) node [block] (IRayleigh) {Inverse Rayleigh}
            (-11,0) node [indic] (5p) {5 parameters}
            (-11,-3) node [indic] (4p) {4 parameters}
            (-11,-6) node [indic] (3p) {3 parameters}
            (-11,-9) node [indic] (2p) {2 parameters}
            (-11,-12) node [indic] (1p) {1 parameter};                            
  \draw[dashed] (IF) -- node[text width=4em, text centered]{$p=0$}(IF1);
  \draw[dashed] (IF) -- node[text width=4em, text centered]{$b=1$}(F);
  \draw[dashed] (IF) -- node[text width=4em, text centered]{$p\to \infty$}(IF2);
  \draw[dashed] (IF1) -- node[text width=4em, text centered]{$q=1$}(Par3);
  \draw[dashed] (IF1) -- node[text width=4em, text centered]{$x_0=0$}(TB12);
  \draw[dashed] (IF1) -- node[text width=4em, text centered]{$b=1$}(Par2);
  \draw[dashed] (F) -- node[text width=4em, text centered]{$p=0$}(Par2);
  \draw[dashed] (F) -- node[text width=4em, text centered]{$x_0=0$}(GL);
  \draw[dashed] (F) -- node[near end, text width=8em, text centered]{$x_0=c(p+1)^{-\frac{1}{q}}$}(Stoppa);
  \draw[dashed] (F) -- node[text width=4em, text centered]{$p\to \infty$}(Frechet); 
  \draw[dashed] (IF2) -- node[text width=4em, text centered]{$b=1$}(Frechet);
  \draw[dashed] (Par3) -- node[text width=4em, text centered]{$x_0=0$}(Fisk);
  \draw[dashed] (TB12) -- node[text width=4em, text centered]{$q=1$}(Fisk);
  \draw[dashed] (TB12) -- node[text width=4em, text centered]{$c=1$}(Burr3);
  \draw[dashed] (TB12) -- node[text width=4em, text centered]{$b=1$}(Lomax); 
  \draw[dashed] (Par2) -- node[text width=4em, text centered]{$x_0=0$}(Lomax); 
  \draw[dashed] (Par2) -- node[near end, text width=4em, text centered]{$x_0=c$}(Par1);
  \draw[dashed] (GL) -- node[near end, text width=4em, text centered]{$p=0$}(Lomax); 
  \draw[dashed] (Stoppa) -- node[near end, text width=4em, text centered]{$p=0$}(Par1);
  \draw[dashed] (GL) -- node[near end, text width=4em, text centered]{ $p \to \infty$}(Gumbel); 
  \draw[dashed] (Stoppa) -- node[text width=4em, text centered]{$p \to \infty$}(Gumbel);
  \draw[dashed] (Frechet) -- node[text width=4em, text centered]{$x_0=0$}(Gumbel); 
  \draw[dashed] (Gumbel) -- node[text width=4em, text centered]{$q=1$}(IExponential);
  \draw[dashed] (Gumbel) -- node[text width=4em, text centered]{$q=2$}(IRayleigh);                   
\end{tikzpicture}
\end{sideways}
\end{center}

\newpage


\section{Main properties} \label{Properties}

In this section we present and discuss the main properties of the IF distribution. Since it contains so many special cases, the subsequent results provide in a single sweep those properties for the various distributions mentioned in Section~\ref{sec:special}.


\subsection{Cumulative distribution function, survival function and hazard function}
One major advantage of the IF distribution is that the cumulative distribution function (cdf) can be written under closed form:
\begin{align}
	F_p(x)
	 =\left \{ 
	 	\begin{array}{l l} 
			\left(1-\frac{1}{p+1}G_p(x)^{-q}\right)^{p+1} \quad &\text{if } b>0, \\
			1-\left(1-\frac{1}{p+1}G_p(x)^{-q}\right)^{p+1} \quad &\text{if } b<0.  \label{CDF}
		\end{array}
		\right.
\end{align}
Consequently, the survival or reliability function $S_p(x)=1-F_p(x)$ is extremely simple, too. The same holds true for the hazard function, defined as the quotient of the pdf and the survival function:
\begin{align*}
	H_p(x)=\frac{f_p(x)}{S_p(x)}
	= \left \{ 
	 	\begin{array}{l l} 
			q~g_p(x)G_p(x)^{-q-1}\frac{\left(1-\frac{1}{p+1}G_p(x)^{-q}\right)^{p}}{1-\left(1-\frac{1}{p+1}G_p(x)^{-q}\right)^{p+1}} \quad &\text{if } b>0, \\
			-q~g_p(x)G_p(x)^{-q-1} \frac{1}{\left(1-\frac{1}{p+1}G_p(x)^{-q}\right)}\quad &\text{if } b<0 .
		\end{array}
		\right.
\end{align*}


\subsection{Quantile function and median}

Very conveniently, the quantile function takes a nice form thanks to the simple expression of the cdf~\eqref{CDF}. Given the wide range of quantile-based statistical tools and methods such as QQ-plots, interquartile range or quantile regression, this is a very welcomed feature of the IF distribution.
For $b>0$, the quantile function is given by
\begin{align*}
	Q_p^+(y)=F_p^{-1}(y) = \left\{
				\begin{array}{l l}
					x_0+c\left((1-y)^{-\frac{1}{q}}-1\right)^{\frac{1}{b}} &\text{if } p=0,\\
					x_0+c(p+1)^{-\frac{1}{bq}}\left(\left(1-y^{\frac{1}{p+1}}\right)^{-\frac{1}{q}}-1\right)^{\frac{1}{b}} &\text{if }  0< p<\infty, \\
					x_0+c\left(\ln \left(\frac{1}{y}\right)\right)^{-\frac{1}{bq}} &\text{if } p \to \infty,
				\end{array}
			\right.
\end{align*}
for $y\in[0,1]$. The expression for $b<0$ is readily obtained via the relationship $Q_p^-(y)=Q_p^+(1-y)$, and we define the quantile function $Q_p(y)$ as $Q_p^+(y)$ if $b>0$ and as $Q_p^-(y)$ if $b<0$. The median is  uniquely defined as 
\begin{align*}
	\mbox{Median}= \left\{
				\begin{array}{l l}
					x_0+c\left((2^{\frac{1}{q}}-1\right)^{\frac{1}{b}} &\text{if } p=0,\\
					x_0+c(p+1)^{-\frac{1}{bq}}\left(\left(1-2^{-\frac{1}{p+1}}\right)^{-\frac{1}{q}}-1\right)^{\frac{1}{b}} &\text{if }  0< p<\infty, \\
					x_0+c \left( \ln \left(2\right)\right)^{-\frac{1}{bq}} &\text{if } p \to \infty.
				\end{array}
			\right.
\end{align*}


\subsection{Random variable generation}
The closed form of the quantile functions entails a straightforward  random variable generation process from the IF. Indeed, it suffices to generate a random variable $U$ from a uniform distribution on the interval~$[0,1]$, and then apply $Q_p$ to it. The resulting random variable $Q_p(U)$ follows the IF~distribution. The simplicity of the procedure is particularly important for Monte Carlo simulation purposes.


\subsection{Moments}

Mean and variance of the IF1, IF2 and IF3 distributions can be written out explicitly in terms of the Gamma and Beta functions. 
\begin{itemize}
	\item Mean and variance of the IF1 distribution ($p=0$) are given by
	\begin{equation}\label{mean}\mathbb{E}(X)=x_0+cq~B\left(q-\frac{1}{b},1+\frac{1}{b}\right) \quad \text{if } \left \{ \begin{array}{ll}
																	b>0 \text{ and } 1<bq, \\
																	b<0 \text{ and } 1<-b,\\
																	\end{array}
																	\right.\end{equation}
	and																
	$$\mathbb{V}(X)=c^2\left[q~B\left(q-\frac{2}{b},1+\frac{2}{b}\right)-\left(q~B\left(q-\frac{1}{b},1+\frac{1}{b}\right)\right)^2\right]\quad \text{if } \left \{ \begin{array}{ll}
																	b>0 \text{ and } 2<bq, \\
																	b<0 \text{ and } 2<-b.\\
																	\end{array}
																	\right.$$															
	\item Mean and variance of the IF2 distribution ($p\to\infty$) are given by
	$$\mathbb{E}(X)=x_0+c~\Gamma\left(1-\frac{1}{bq}\right) \quad \text{if } \left \{ \begin{array}{ll}
																	b>0 \text{ and } 1<bq, \\
																	b<0,\\
																	\end{array}
																	\right.$$
	and																
	$$\mathbb{V}(X)=c^2 \left[\Gamma\left(1-\frac{2}{bq}\right)-\left(\Gamma\left(1-\frac{1}{bq}\right)\right)^2\right] \quad \text{if } \left \{ \begin{array}{ll}
																	b>0 \text{ and } 2<bq, \\
																	b<0.&\\
																	\end{array}
																	\right.$$																	
	\item Mean and variance of the IF3 distribution ($0<p< \infty$ and $b=1$) are
	$$\mathbb{E}(X)=x_0+c (p+1)^{1-\frac{1}{q}} \left(B\left(1-\frac{1}{q},p+1\right)-\frac{1}{p+1}\right) \quad \text{if } 1<q,$$
	and
	\begin{eqnarray*}
	\mathbb{V}(X)&=&c^2(p+1)^{1-\frac{2}{q}}\left[B\left(1-\frac{2}{q},p+1\right)-\frac{1}{p+1}\right]\\
	&&+c^2(p+1)^{2-\frac{2}{q}}\left[B\left(1-\frac{1}{q},p+1\right)-\frac{1}{p+1}\right]^2 \quad \text{if } 2<q.		\end{eqnarray*}	\end{itemize}

To get a flavor of the underlying calculations, we stress that, e.g., expression~\eqref{mean}  is best obtained by first performing the change of variable $y= \left(1+\left(\frac{x-x_0}{c}\right)^b\right)^{-1}$. For higher-order  moment expressions, we refer the interested reader to~\cite{sinner2016}, available on the ArXiv. That paper also discusses  entropic properties of the Interpolating Family.


 \subsection{Modality}\label{modcalc}

Determining the mode of a distribution is an important issue, which we tackle in this section. The detailed calculations are deferred to the Appendix. We study the derivative of $x\mapsto f_p(x)$, with particular emphasis on the three main subfamilies IF1, IF2 and~IF3 described in Section~\ref{sec:special}.

\hfill

The derivative of the pdf vanishes either at the boundary $x=x_0$ of the domain or at 
\begin{equation*}
x= x_0+c (p+1)^{-\frac{1}{bq}} \left(t^{-\frac{1}{q}}-1\right)^{\frac{1}{b}},
\end{equation*}
where $t$ is solution of the almost cyclic equation
\begin{equation}\label{cyclicequation}
(b-1) t^{-\frac{1}{q}} (1-t) -b(q+1)(t^{-\frac{1}{q}}-1)(1-t) + pbq (t^{-\frac{1}{q}}-1) t =0.
\end{equation}
This allows us to draw the following conclusions regarding the modality of the IF~distribution.
\begin{itemize} 
\item The mode of the IF1 distribution ($p=0$) is given by
\begin{equation*}
\left\{
\begin{array}{ll}
x_0 & \text{ if }b = -\frac{1}{q} \text{ or } b=1,\\
x_0+c\left(\frac{b-1}{bq+1}\right)^{\frac{1}{b}} & \text{ if } b<-\frac{1}{q} \text{ or } b>1,
\end{array}
\right.
\end{equation*}
whereas in the remaining cases, i.e.~$b \in ]-\frac{1}{q};1[$, there is a vertical asymptote in $x=x_0$. We plot in Figure~\ref{fig1} a contour plot of the mode of the IF1~distribution.
\begin{figure}[h!]
	\begin{center}
	\includegraphics[scale=0.6]{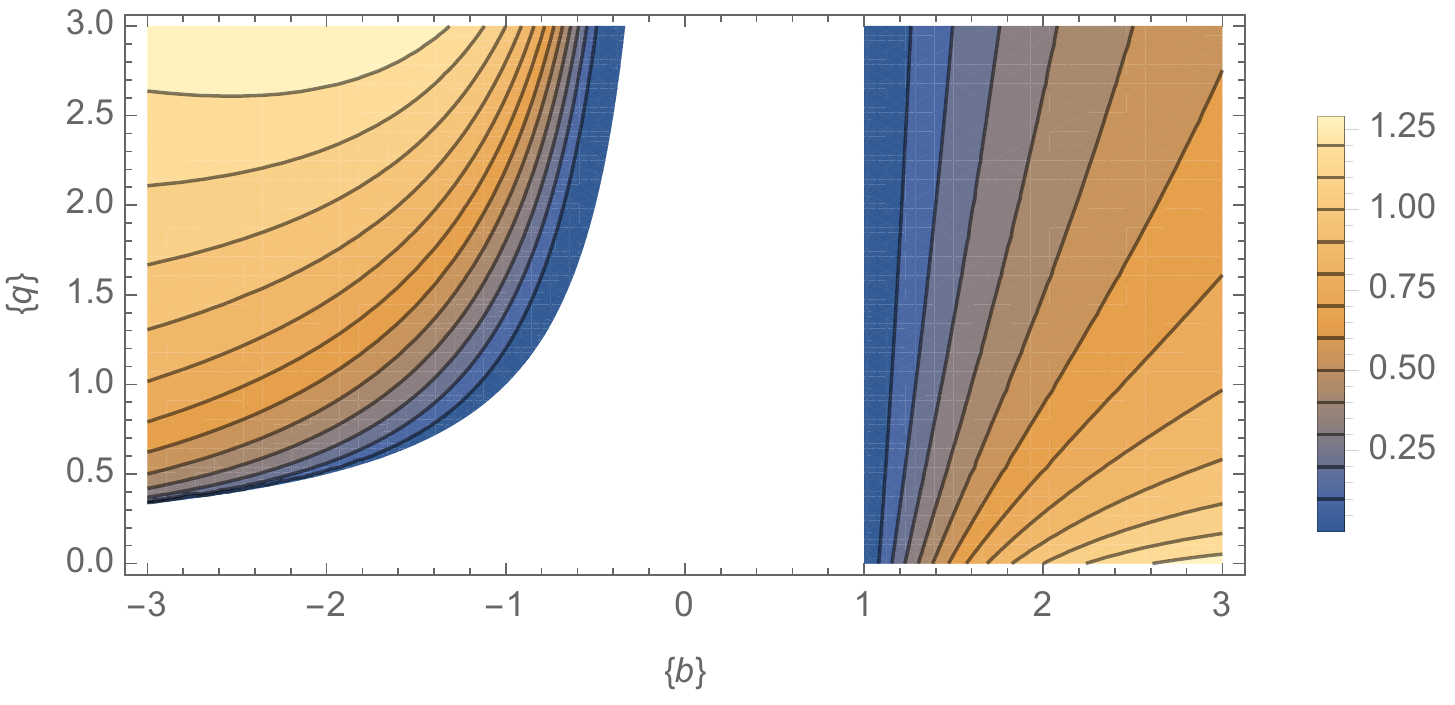}
	\caption{Contour plot of the mode of the IF1 distribution with parameters $c=1$ and $x_0=0$.}\label{fig1}
	\end{center}  
\end{figure}
\item The mode of the IF2 distribution ($p\to \infty$) is given by
\begin{equation*}
\left\{
\begin{array}{ll}
x_0 & \text{ if }b = -\frac{1}{q},\\
x_0+c\left(\frac{bq}{bq+1}\right)^{\frac{1}{bq}} & \text{ if } b<-\frac{1}{q} \text{ or } b>0,
\end{array}
\right.
\end{equation*}
whereas in the remaining cases, i.e.~$b \in ]-\frac{1}{q},0]$, there is a vertical asymptote in $x=x_0$. Figure~\ref{fig2} shows a contour plot of the mode of the IF2 distribution.
\begin{figure}[h!]
	\begin{center}
	\includegraphics[scale=0.6]{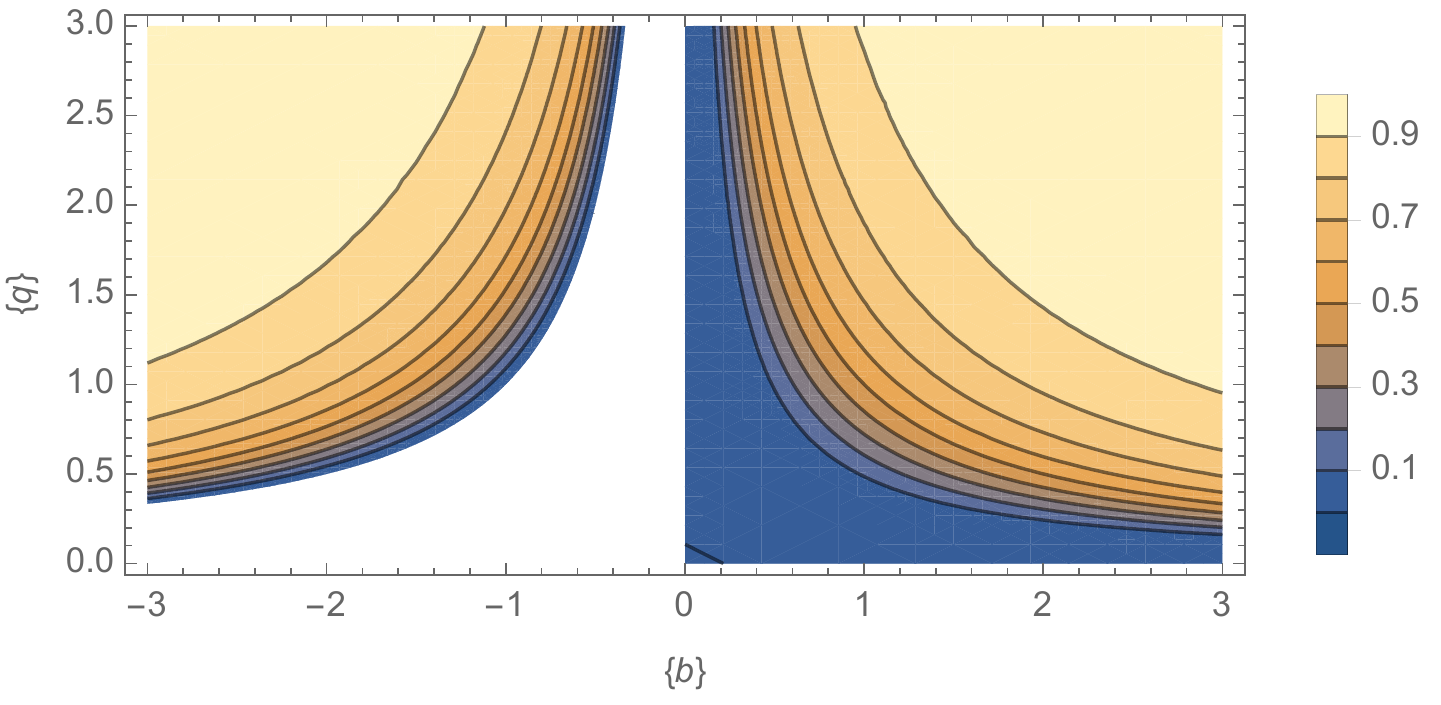}
	\caption{Contour plot of the mode of the IF2 distribution with parameters $c=1$ and $x_0=0$.}\label{fig2}
	\end{center}  
\end{figure}
\item The mode of the IF3 distribution ($0<p<\infty$ and $b=1$) is given by
$$
x_0+c(p+1)^{-\frac{1}{q}}\left(\left(\frac{q+1}{(p+1)q+1}\right)^{-\frac{1}{q}}-1\right).
$$
A contour plot of the mode of the IF3 distribution can be seen in Figure~\ref{fig3}.
\begin{figure}[h!]
	\begin{center}
	\includegraphics[scale=0.35]{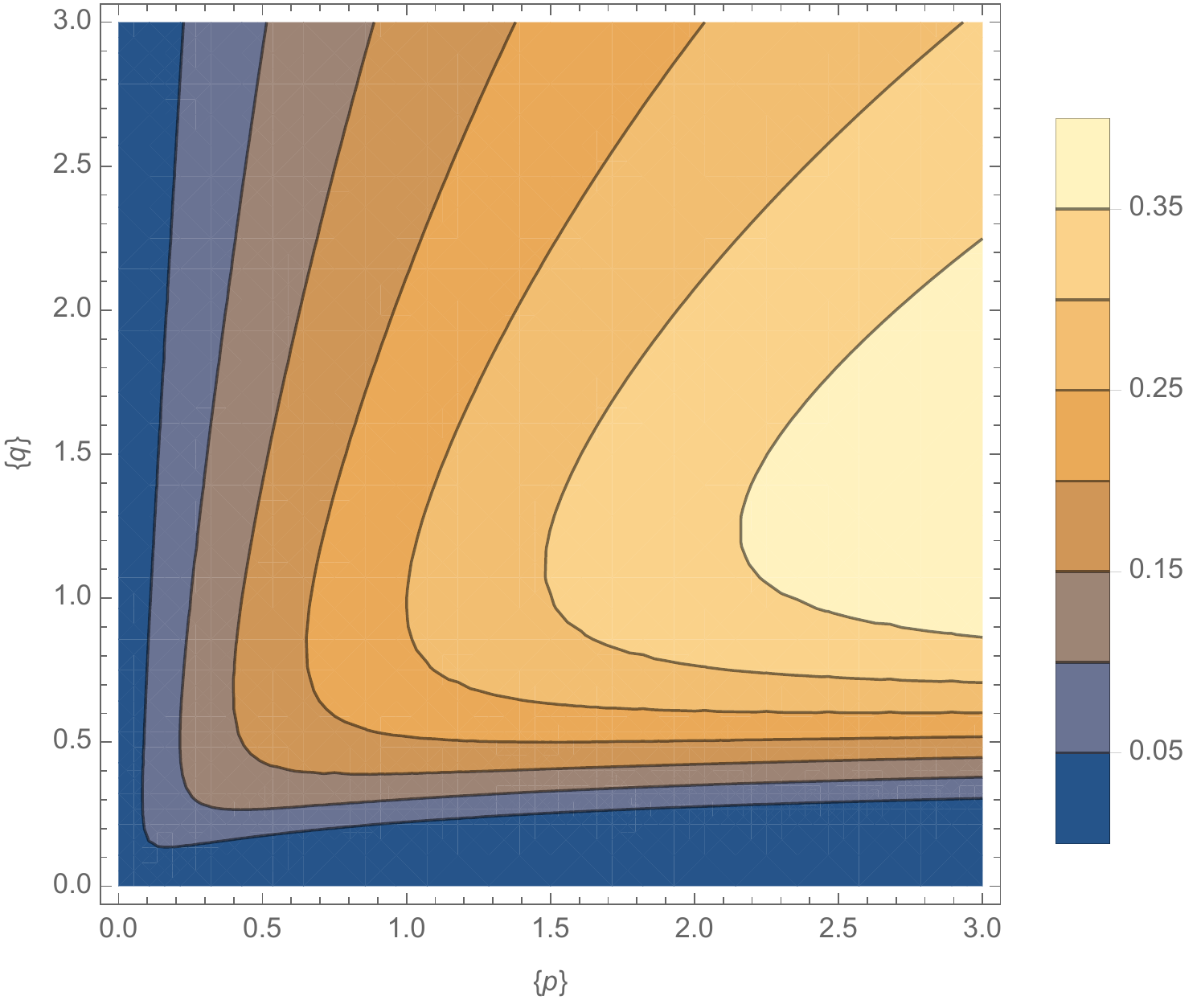}
	\caption{Contour plot of the mode of the IF3 distribution with parameters $c=1$ and $x_0=0$.}\label{fig3}
	\end{center}  
\end{figure}
\end{itemize}

We see that all three subfamilies are unimodal, which is coherent with the related special cases from the literature. Moreover, we have derived the exact expressions of the modes. This unimodality is a very attractive feature from an interpretability point of view: bimodal or multimodal distributions are arguably best modeled as a mixture of unimodal distributions. It is therefore not surprising that many new distributions are built with the target of being unimodal; for recent examples, see e.g.~\cite{jones2014generating, katojones2015, fujisawaabe2015}.


\section{Inferential aspects}


\subsection{Parameter estimation}\label{MLEclassic}

Let $(x_1,x_2,\dots,x_n)$ be  independent and identically distributed  observations from a population with density \eqref{IF}. Then the log-likelihood function can be expressed as
\begin{eqnarray*}
	\ell\left(p,b,c,q,x_0\right)&=&  n \ln \left(\frac{|b|q}{c^b} \right)+(b-1) \sum \limits_{i=1}^n \ln(x_i-x_0) \\ &&-(q+1) \sum \limits_{i=1}^n\ln  \left((p+1)^{-\frac{1}{q}} + \left(\frac{x_i-x_0}{c}\right)^b \right) \\ &&+p  \sum \limits_{i=1}^n \ln \left(1-\frac{1}{p+1}\left((p+1)^{-\frac{1}{q}} + \left(\frac{x_i-x_0}{c}\right)^b \right)^{-q}\right).
\end{eqnarray*}
Deriving the score functions is straightforward by differentiating the log-likelihood function with respect to each of the parameters $p,b,c,q,x_0$. To estimate the parameter $x_0$ we impose the constraint that it is smaller than or equal to the lowest value of the data set. 
The related likelihood equations can readily be solved by any numerical root-finding algorithm. We carried out the calculations in Section \ref{sec:Applications} with \textsl{Wolfram Mathematica 10}. More precisely, we used the function \textsl{NMaximize} with the numerical maximization algorithm \textsl{Random search method} enhanced with the option \textsl{InteriorPoint} and we limited the maximum number of iterations to $10^6$.


\subsection{Submodel testing}\label{submod}

We have shown in Section \ref{sec:special} the many distinct submodels that the IF  nests. Consequently, it is natural to propose tests for submodels within the IF. This can be done by likelihood ratio tests. For each parameter $\eta\in\{p,b,c,q,x_0\}$, we denote by $\hat\eta$ the unconstrained maximum likelihood estimate
and by $\hat\eta_r$ the maximum likelihood estimate under the restricted submodel of interest. For example, testing for the Pareto type~I distribution against the larger IF~model can be achieved by the test statistic
$$
T_{\rm Pareto}=-2\left(\ell(0,1,\hat c_r,\hat q_r,\hat c_r)-\ell(\hat p,\hat b,\hat c,\hat q,\hat x_0)\right)
$$
rejecting $\mathcal{H}_0:\{b=1\}\cap\{x_0=c\}$ at asymptotic level $\alpha$ whenever $T_{\rm Pareto}$ exceeds~$\chi^2_{3;1-\alpha}$, the $\alpha$-upper quantile of the chi-squared distribution with three degrees of freedom. Similarly, testing for the Weibull against the IF distribution leads to the test statistic
$$
T_{\rm Weibull}=-2\left(\ell(\infty,-1,\hat c_r,\hat q_r,0)-\ell(\hat p,\hat b,\hat c,\hat q,\hat x_0)\right)
$$
to be compared with $\chi^2_{3;1-\alpha}$. Similar submodel tests, against other families in which the Pareto and Weibull are respectively nested, are proposed in~\cite{falk2008} and \cite{mudholkar1996generalization}, to cite but these.


\subsection{Dealing with censored data}\label{censor}

In survival analysis, the data are most often sampled from a population containing censored observations. 
Depending on when this phenomenon occurs, the observations can be attributed to one of the four following sets: Obs=\{the phenomenon occurs during the study and the observation is uncensored\}, Left=\{the phenomenon occurs before the start of the study\}, Int=\{the phenomenon occurs within a finite interval of time\}, or Right=\{the phenomenon does not occur during the period of the study\}. 
The contribution of censored observations to the likelihood function is given by
\begin{eqnarray*}
	\mathcal{L}\left(p,b,c,q,x_0\right)&=& \prod \limits_{j \in \text{Obs}} f_{p}(x_j) \prod \limits_{j \in \text{Left}} F_{p}(x_j) \prod \limits_{j \in \text{Right}} \left(1-F_{p}(x_j)\right) \\ &&\prod \limits_{j \in \text{Int}} \left(F_{p}(x_j^R)-F_{p}(x_j^L)\right).
\end{eqnarray*}
The maximum likelihood estimates for the parameters are then derived in a similar fashion as for the non-censored data of Section~\ref{MLEclassic}. The very simple closed-form expression of the cdf, see Section \ref{Properties}, seems tailor-made for dealing with censored data, especially in comparison to the Generalized Beta family where the cdf needs to be computed numerically.


\section{Applications}\label{sec:Applications}

We illustrate the new family of size distributions by three applications. In the first two applications, one on liability claims data and one on fire size data, we compare the fitting capacities of the Interpolating Family of size distributions to the Generalized Beta family. Our means of comparison shall be the Akaike Information Criterion (AIC) and the Bayesian Information Criterion (BIC). Given that the Pareto and Weibull distributions are natural choices to model respectively the first and second data set, we shall  test for these submodels within the IF model by means of the tests of Section~\ref{submod}. The third application on survival data shall illustrate how well the IF is suited for dealing with censored data.


\subsection{Actuarial science}

\begin{center}
	\begin{figure}[ht]
		\includegraphics[scale=0.55]{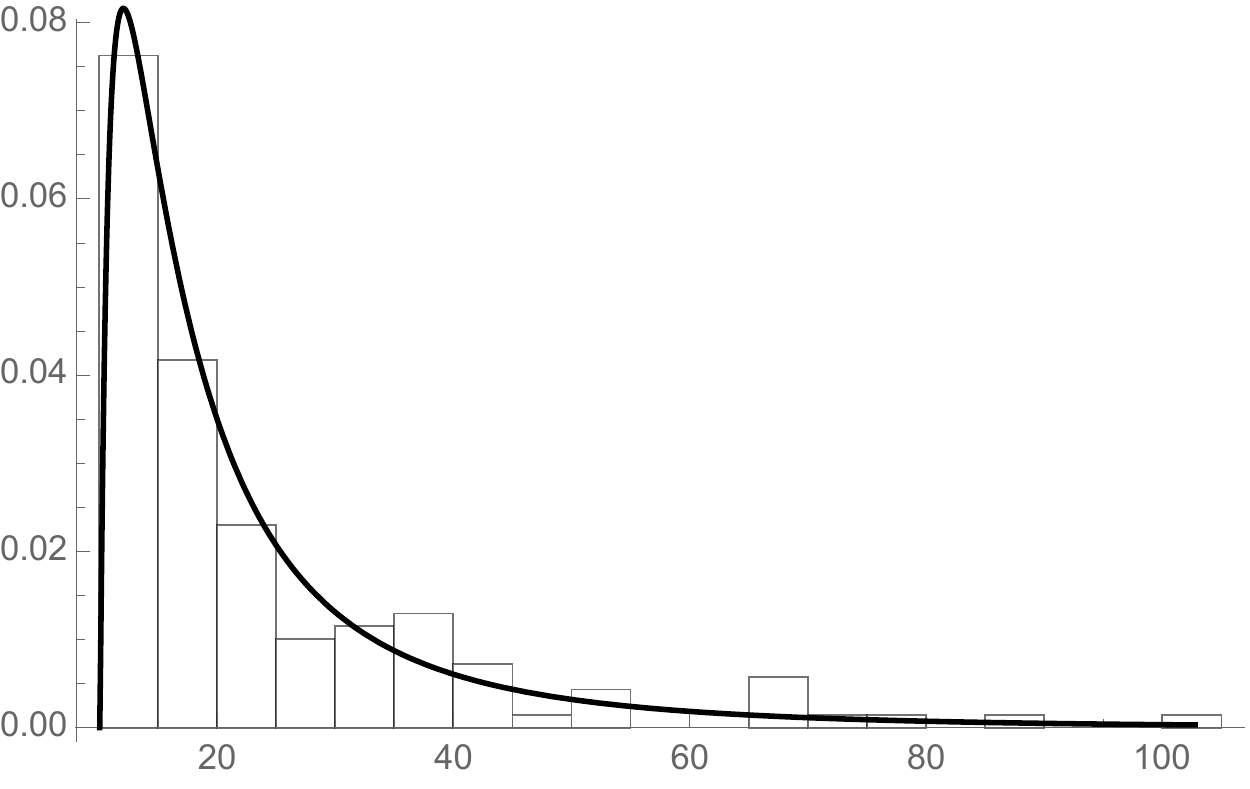}
		\caption{Histogram of the actuarial liability claims data (in Euro, and divided by $10^5$) together with the best-fitting IF density.}
		\label{ACTU}
	\end{figure}
\end{center}

We will analyze $n=139$ motor third party liability claims, recorded during the period 1998--2013. These data are provided by a European reinsurance company. Their values can be appreciated from Figure~\ref{ACTU}.

Table \ref{AICBIC} shows the results of the maximum log-likelihood (MLL) values as well as the AIC and BIC values obtained from all the models used in this analysis, while Tables \ref{ParametersIF} and \ref{ParametersGB} contain the corresponding parameter estimates.
The IF distribution has lower AIC and BIC values than the GB distribution and thus fits best the data. This fit is illustrated on Figure \ref{ACTU}.

 \hfill 

\begin{center}
	\begin{tabular}{|c|c|c|c|c|}
	\hline
	Distribution & $k$ & MLL & AIC & BIC\\
	\hline
	\hline
	IF  & 5 & -499.255 & 1006.510 & 1018.248 \\
	\hline
	\hline
	GB  & 5 & -499.671 & 1009.342 & 1024.014 \\
	\hline
	\hline
	IF1  & 4 & -500.200 & 1008.400 & 1020.138 \\
	\hline
	IF2  & 4 & -502.034 & 1012.068 & 1023.806 \\
	\hline
	IF3  & 4 & -500.069 & 1008.138 & 1019.876 \\
	\hline
	\hline
	GB1  & 4 & -500.210 & 1008.420 &1020.158 \\
	\hline
	GB2  & 4 & -500.186 & 1008.372 & 1020.110\\
	\hline	
	\end{tabular}
\captionof{table}{\label{AICBIC} Maximum value of the log-likelihood function as well as the AIC and BIC values for the actuarial liability claims data set.}
\end{center}

\begin{center}
	\begin{tabular}{|c|c|c|c|c|c|}
	\hline
	Distribution & $\hat{p}$ & $\hat{b}$ & $\hat{c}$ & $\hat{q}$ & $\hat{x}_0$\\
	\hline
	\hline
	IF  & 28.593 & -0.281 & 0.004 & 30.720 & 10.059 \\
	\hline
	IF1  &  0 & 23.368 & 11.447 & 0.073 & 0\\
	\hline
	IF2  & $\infty$ & 0.310 & 7.108 & 4.469 & 8.106  \\
	\hline
	IF3  &  0.261 & 1 & 25.199 & 2.931 & 10.250 \\	
	\hline	
	\end{tabular}
\captionof{table}{\label{ParametersIF}  Parameter estimation values of the IF distribution and of its submodels for the actuarial liability claims data.}
\end{center}

\begin{center}
	\begin{tabular}{|c|c|c|c|c|c|}
	\hline
	Distribution & $\hat{a}$ & $\hat{b}$ & $\hat{c}$ & $\hat{p}$ & $\hat{q}$\\
	\hline
	\hline
	GB  & -0.421 & 515.021 & 0.808 & 19.166 & 1.307\\
	\hline
	GB1  &  -0.034 & 10.164 & 0 & 68.525 & 1.651\\
	\hline
	GB2  &  -27.680 & 11.652 & 1 & 0.062 & 0.733\\
	\hline	
	\end{tabular}
\captionof{table}{\label{ParametersGB}  Parameter estimation values of the GB distribution and of its submodels for the actuarial liability claims data.}
\end{center}

%
%

A popular choice to model this kind of data is the Pareto distribution, see for example~\cite{benktander1970schadenverteilung}. The MLL value of the Pareto type~I model is $-506.417$, and the likelihood ratio test statistic for this distribution with respect to the IF~distribution equals $T_{\text{Pareto}}=-2(-506.417+499.255)=14.324$, yielding a $p$-value of~$0.0025$. Consequently we clearly reject the Pareto submodel in favor of the IF~distribution.
 

\subsection{Environmental science}

The second data set consists of $n=102$ means of log-burned areas of 26870 wildfires, which were grouped by watersheds. These data were recorded in Portugal during the period 1985--2005. For more information about the data, we refer to~\cite{barros2012identifying} and~\cite{garcia2014test, garcia2014central}. 

\hfill

Table \ref{AICBICFire} shows the results of the maximum log-likelihood (MLL) values as well as the AIC and BIC values obtained from all the models used in this analysis, while Tables \ref{ParametersIFFire} and \ref{ParametersGBFire} contain the corresponding parameter estimates. As we can see, the IF1 distribution happens to coincide with the best possible IF model, hence has lowest AIC and BIC values within the Interpolating Family. The IF1 also outperforms the best Generalized Beta model, namely the GB2.  

\hfill

\begin{center}
	\begin{tabular}{|c|c|c|c|c|}
	\hline
	Distribution & $k$ & MLL & AIC & BIC\\
	\hline
	\hline
	IF & 5 &  -575.591 & 1161.182 & 1174.307\\
	\hline
	\hline
	GB  & 5 &  -576.066 & 1162.132 & 1175.257\\
	\hline
	\hline
	IF1  & 4 & -575.591 & 1159.182 & 1169.682 \\
	\hline
	IF2  & 4 &  -577.305 & 1162.610 & 1173.110\\
	\hline
	IF3  & 4 &  -576.051 & 1160.102 & 1170.602\\
	\hline
	\hline
	GB1  & 4 &  -581.418 & 1170.836 & 1181.336\\
	\hline
	GB2  & 4 & -576.270 & 1160.540 & 1171.040\\
	\hline	
	\end{tabular}
\captionof{table}{\label{AICBICFire} Maximum value of the log-likelihood function as well as the AIC and BIC values for the environmental data.}
\end{center}

\begin{center}
	\begin{tabular}{|c|c|c|c|c|c|}
	\hline
	Distribution & $\hat{p}$ & $\hat{b}$ & $\hat{c}$ & $\hat{q}$ & $\hat{x}_0$\\
	\hline
	\hline
	IF  &  0 & -1.769 & 76.503 & 0.801 & 15.357\\
	\hline
	IF1  &  0 & -1.769 & 76.503 & 0.801 & 15.357\\
	\hline
	IF2  & $\infty$ &   0.845 & 57.814 & 1.662 & 0\\
	\hline
	IF3  & 1.379  & 1 & 109.968 & 2.150 & 12.996  \\	
	\hline	
	\end{tabular}
\captionof{table}{\label{ParametersIFFire} Parameter estimation values of the IF distribution and of its submodels for the environmental data.}
\end{center}

\begin{center}
	\begin{tabular}{|c|c|c|c|c|c|}
	\hline
	Distribution & $\hat{a}$ & $\hat{b}$ & $\hat{c}$ & $\hat{p}$ & $\hat{q}$\\
	\hline
	\hline
	GB  & -1.393 & 41.496 & 0.834 & 1.259 & 2.436 \\
	\hline
	GB1  &  0.059 & $8.791 \cdot 10^8$ & 0 & 248.824 & 392.153 \\
	\hline
	GB2  &  0.875 & 0.416 & 1 & 198.305 & 2.355\\
	\hline	
	\end{tabular}
\captionof{table}{\label{ParametersGBFire} Parameter estimation values of the GB distribution and its submodels for the environmental data.}
\end{center}

The classical choice to model environmental data  is the Weibull distribution.
The likelihood ratio test statistic for the Weibull distribution with respect to the IF~distribution yields $T_{\text{Weibull}}=-2(-598.87+575.591)=46.558$  and a $p$-value of~$0.00$, hence we reject the Weibull submodel in favor of the broader IF model. 


\subsection{Survival analysis}

For the previous two data sets we have seen that our IF distribution outperforms the GB distribution as well as the classical choices, the Pareto Type I and Weibull distributions. We shall now investigate its fitting abilities on censored data, for which we have already argued in Section~\ref{censor} that the IF presents important computational advantages over the GB distribution. Thus we only focus on the IF here and compare it to the Weibull, a popular choice in survival analysis.

\hfill

The third data set consists of $n=228$ survival data from the North Central Cancer Treatment Group. They measure the survival time, in days, of patients suffering from advanced lung cancer. These data are available in the `survival'-package in R. The phenomenon of interest is the death of the patients and we are consequently facing  right-censored data whenever a patient has survived longer than the allocated study time. This results in $n_{\text{Obs}}=165$ and $n_{\text{Right}}=63$ data.  
The maximum log-likelihood value of the IF model is $-1153.460$, whilst the one of the  Weibull model is $-1153.850$. This yields a likelihood ratio test value of $T_{\text{Weibull}}=-2(-1153.850+1153.460)=0.78$ and a $p$-value of $0.854$. Thus in this case we do not reject the null hypothesis. In what follows, we shall therefore further compare the outcomes of the IF to those of the Weibull.

\hfill

In survival analysis, a question of interest is to study the survival time of patients suffering from lung cancer. 
We computed the probability that a patient survives more than $t=1,2,3$ and $4$ years (with a daycount of 365) with both the IF~model and the Weibull model. The results are given in Table \ref{Survival} below. 

\hfill

\begin{center}
	\begin{tabular}{|c|c|c|}
	\hline
	$t$ & IF & Weibull \\
	\hline
	\hline
	1 & 0.4384 & 0.4330 \\
	\hline
	2 & 0.1256 & 0.1243 \\
	\hline
	3 & 0.0218 & 0.0285 \\
	\hline
	4& 0.0015 & 0.0055 \\
	\hline	
	\end{tabular}
\captionof{table}{\label{Survival} Probability to survive $t=1,2,3$ and $4$ years on basis of the lung cancer data.}
\end{center}

A second quantity of interest is the probability to decease before time $t+1$ knowing that the patient was alive at time $t$. This can be computed as $$\mathbb{P}(X<t+1~|~X>t)=\frac{F(t+1)-F(t)}{1-F(t)}$$ and is again readily computable from the IF. Table \ref{Proba} shows the results. The Weibull model yields lower probabilities to decease within one year starting from time $t=2$. 

\hfill

\begin{center}
	\begin{tabular}{|c|c|c|}
	\hline
	& IF & Weibull \\
	\hline
	\hline
	 $\mathbb{P}(X<1~|~X>0)$ &0.5615&0.5670\\
	\hline
	$\mathbb{P}(X<2~|~X>1)$ &0.7134&0.7130 \\
	\hline
	 $\mathbb{P}(X<3~|~X>2)$ &0.8262&0.7704\\
	\hline
	 $\mathbb{P}(X<4~|~X>3)$ &0.9292&0.8057\\
	\hline	
	$\mathbb{P}(X<5~|~X>4)$ &0.9934&0.8304 \\
	\hline
	\end{tabular}
\captionof{table}{\label{Proba} Probability to decease before time $t+1$ knowing that the patient was alive at time $t$ with respect to the lung cancer  data.}
\end{center}

This third data set shows that, also in survival analysis, opting for the IF distribution is a  good choice. The Weibull distribution is more parameter-parsimonious and cannot be significantly rejected in favor of the IF, nevertheless the latter has the higher MLL value and its probabilities  from Tables~\ref{Survival} and~\ref{Proba} yield  precious information about the data under investigation. In general, it is safer to bet on the ~IF as it yields an excellent fit for very diverse data sets.


\section{Conclusion}

In this paper we introduced the Interpolating Family of size distributions, and studied its stochastic and inferential properties. Three very distinct real data sets allowed to appreciate its excellent fitting abilities, also compared to the Generalized Beta distribution. We mentioned in the Introduction that this new distribution should be a viable alternative to the GB distribution. To help the reader choose his/her favorite distribution, we conclude the paper with a short comparison between the existing Generalized Beta distribution and the new Interpolating Family. A common feature is the high flexibility as they depend on five parameters and they both nest power laws as well as power laws with exponential cut-off as special cases. The GB distribution has the advantage to nest interesting special cases not covered by the IF distribution as, for instance, Gamma-type distributions. Major advantages of the IF distribution over the GB~distribution are the closed-form cumulative distribution function, very tractable quantile expressions and a simple random variable generation process. Consequently, the IF~distribution can readily be applied to censored data in survival analysis. Moreover, the normalizing constant does not involve any special function. Finally, the construction of the IF distribution is very geometric as it arises naturally when interpolating between power laws and power laws with exponential cut-off.

\section*{Appendix: Mode calculation}

In this Appendix we detail some of the steps we skipped in Section~\ref{modcalc} to calculate the mode of the Interpolating Family of size distributions. 
The derivative of the pdf~(\ref{IF}) vanishes if and only if \small
\begin{eqnarray*}
	0=&(b-1) \left(\frac{x-x_0}{c}\right)^{b-2} \left((p+1)^{-\frac{1}{q}}+\left(\frac{x-x_0}{c}\right)^b\right)^{-q-1} \left(1-\frac{1}{p+1}\left( (p+1)^{-\frac{1}{q}}+\left(\frac{x-x_0}{c}\right)^b\right)^{-q}\right)^p \\
	&-b(q+1)\left(\frac{x-x_0}{c}\right)^{2b-2} \left((p+1)^{-\frac{1}{q}}+\left(\frac{x-x_0}{c}\right)^b\right)^{-q-2} \left(1-\frac{1}{p+1}\left( (p+1)^{-\frac{1}{q}}+\left(\frac{x-x_0}{c}\right)^b\right)^{-q}\right)^p \\
	&+\frac{pbq}{p+1} \left(\frac{x-x_0}{c}\right)^{2b-2} \left((p+1)^{-\frac{1}{q}}+\left(\frac{x-x_0}{c}\right)^b\right)^{-2q-2}\left(1-\frac{1}{p+1}\left( (p+1)^{-\frac{1}{q}}+\left(\frac{x-x_0}{c}\right)^b\right)^{-q}\right)^{p-1}. 
\end{eqnarray*} \normalsize
If we set $y= \frac{x-x_0}{c}$, then the above holds true if either 
\small \begin{equation*}
y^{b-2}
\left((p+1)^{-\frac{1}{q}}+y^b\right)^{-q-2} 
\left(1-\frac{1}{p+1}\left( (p+1)^{-\frac{1}{q}}+y^b\right)^{-q}\right)^{p-1}=0
\end{equation*} \normalsize
or
\small \begin{eqnarray}\label{secondequation}
&&0 = (b-1) \left((p+1)^{-\frac{1}{q}}+y^b\right) \left(1-\frac{1}{p+1}\left( (p+1)^{-\frac{1}{q}}+y^b\right)^{-q}\right) \\
& & \quad\quad -b (q+1) y^{b} \left(1-\frac{1}{p+1}\left( (p+1)^{-\frac{1}{q}}+y^b\right)^{-q}\right)  + \frac{pbq}{p+1} y^{b} \left((p+1)^{-\frac{1}{q}}+y^b\right)^{-q}. \nonumber
\end{eqnarray} \normalsize
The only solution that the first equation can possibly admit is $y=0$. This corresponds to $x=x_0$, i.e.~to the boundary of the domain. On the other hand, equation~(\ref{secondequation}) may admit interior solutions. We separate the analysis of equation~(\ref{secondequation}) in two parts: first the case $p$ finite from which we deduce the mode of the IF1 and IF3 subfamilies and second the case $p \to \infty$ which gives the mode of the IF2~distribution. If $p$ is finite, then we set 
$$t=\frac{1}{p+1}\left((p+1)^{-\frac{1}{q}}+y^b\right)^{-q}$$ and equation~(\ref{secondequation}) simplifies to the almost cyclic equation~(\ref{cyclicequation}):
$$(b-1) t^{-\frac{1}{q}}(1-t)-b(q+1)\left(t^{-\frac{1}{q}}-1\right)(1-t)+pbq\left(t^{-\frac{1}{q}}-1\right)t=0.$$
Solving this equation in all generality is possible numerically but we will restrict ourselves to show how to get closed-form solutions for the two subfamilies IF1 and~IF3.
For the IF1 distribution ($p=0$), equation~(\ref{cyclicequation}) further simplifies to
\begin{align*}
	&(b-1)t^{-\frac{1}{q}} (1-t)-b(q+1)\left(t^{-\frac{1}{q}}-1\right)(1-t) =0. 
\end{align*}
While we recover the boundary solution $x=x_0$ if $b>0$, we also find an interior solution $x=x_0+c \left(\frac{b-1}{bq+1}\right)^{\frac{1}{b}}$ if either $b<-\frac{1}{q}$ or $b>1$. Repeating the procedure with the second derivative of the pdf~(\ref{IF}), a straightforward but tedious calculation shows that the interior solution thus found indeed corresponds to a maximum and that the mode occurs on the boundary $x=x_0$ if either $b=-\frac{1}{q}$ or $b=1$.

For the IF3 distribution ($0<p<\infty$ and $b=1$), equation~(\ref{cyclicequation}) further simplifies to
\begin{align*}
	&-(q+1)\left(t^{-\frac{1}{q}}-1\right)(1-t)+pq\left(t^{-\frac{1}{q}}-1\right)t =0. 
\end{align*}
This equation admits two solutions: the boundary solution $x=x_0$ and the interior solution $x=x_0+c (p+1)^{-\frac{1}{q}} \left(\left(\frac{q+1}{(p+1)q+1}\right)^{-\frac{1}{q}}-1\right)$. One can then check that the latter corresponds to the mode of the IF3~distribution and that this mode, and thus the interior solution, moves towards the boundary as $p$ and $q$ tend to zero.

On the other hand, for the IF2 distribution ($p\to \infty$), equation~(\ref{secondequation}) simplifies to
$$0 = (b-1) y^b - b(q+1)y^b + bq  y^{b-bq}.$$
We deduce that the derivative of the pdf of the IF2 vanishes either at the boundary $x=x_0$ if $b>0$ or at the interior point $x=x_0+c\left(\frac{bq}{bq+1}\right)^{\frac{1}{bq}}$ if either $b<-\frac{1}{q}$ or $b>0$. Similarly as for the IF1, tedious second derivative calculations reveal that the interior solution always corresponds to a maximum and that the mode occurs on the boundary if either $b=-\frac{1}{q}$ or $b=0$.

\section*{Acknowledgments}
Yves Dominicy acknowledges financial support from the Fonds National de la Recherche Scientifique, Communaut\'e Fran\c caise de Belgique, via a Mandat de Charg\'e de Recherche FNRS. Patrick Weber acknowledges financial support via an Aspirant grant from the FNRS.

\bibliographystyle{apalike}
\bibliography{refs}

\end{document}